\documentclass[dvips, 12pt, a4paper]{paper}
\usepackage{amsmath,amssymb,amsthm, dsfont,mathtools}
\usepackage{color}
\definecolor{violet}{rgb}{0.5,0,0.8}

\theoremstyle{plain}
\newtheorem{thm}{Theorem}[section]

\theoremstyle{plain}
\newtheorem{cor}[thm]{Corollary}

\theoremstyle{plain}
\newtheorem{lemma}[thm]{Lemma}

\theoremstyle{plain}

\theoremstyle{plain}
\newtheorem{proposition}[thm]{Proposition}

\theoremstyle{plain}
\newtheorem{definition}[thm]{Definition}

\theoremstyle{plain}

\theoremstyle{plain}
\newtheorem{ass}[thm]{Assumption}

\theoremstyle{plain}
\newtheorem{remark}[thm]{Remark}

\begin{document}

\title{Martingale representations in progressive enlargement by multivariate point processes}

\author{Antonella Calzolari \thanks{Dipartimento di Matematica - Universit\`a di Roma
"Tor Vergata", via della Ricerca Scientifica 1, I 00133 Roma,
Italy }  \and Barbara Torti $^*$} \maketitle
\bigskip
\maketitle
\begin{abstract}
We show that all local martingales with respect to the initially enlarged natural filtration
 of a vector of multivariate point processes can be weakly represented up to the minimum among the explosion times of the components.~We also prove that a strong representation holds if any multivariate point process of the vector has
 almost surely infinite explosion time and discrete mark's space.~Then
we  provide a condition under which the components of the multidimensional local martingale driving the strong representation are pairwise orthogonal.
\end{abstract}

\textbf{AMS 2010} 60G48 \, 60G44 \, 60H05 \,
60G55 \,
60G57
\newpage

\section{Introduction}
Martingale representations is a widely studied topic in stochastic analysis since the seventies.~As well-known it deals, given a filtration
$\mathbb{F}=\left(\mathcal{F}_t\right)_{t\geq 0}$, with the problem of representing all $\mathbb{F}$-local martingales as \textit{(vector) stochastic integral} of
 a predictable process with respect to a fixed (multidimensional) local martingale (\textit{strong representation}), possibly with the addition of the  integral of a
 predictable function with respect to a fixed compensated random measure (\textit{weak representation}).\\
 If the  strong representation holds, then the
 local martingale which represents all the others is said to have the \textit{(strong) predictable representation property with respect to  $\mathbb{F}$} ($\mathbb{F}$-PRP)
 (see e.g.~Definition 1.1 page 362 in \cite{he-wang-yan92}).~Actually, as in \cite{co-nu-sc05} and \cite{davis1}, the definition of PRP can be extended to a countable family of  local martingales.~More precisely  a countable family of  local martingales enjoys the PRP if any local martingale can be written as  linear combination
  of stochastic integrals of
  predictable processes, each  driven by  a fixed element of the family.\\
  If the weak representation holds, when the local martingale
  driving the stochastic integral is the continuous local martingale part of a fixed semimartingale and the random measure is its jump measure, then the
  semimartingale is said to satisfy the \textit{weak predictable representation property with respect to  $\mathbb{F}$} ($\mathbb{F}$-WRP)
  (see e.g.~Definition 13.13 page 368 in \cite{he-wang-yan92}).\\
  When the the PRP holds for a multidimensional local martingale with pairwise orthogonal components (or, analogously,  for a discrete family of local martingales with pairwise orthogonal elements)
  the latter form a set of real local martingales called \textit{$\mathbb{F}$-basis}.~In fact in this case the vector stochastic integral is a \textit{componentwise stochastic integral} and any
  $\mathbb{F}$-local martingale can be expressed as sum of stochastic integrals, one integral for each element of the basis (see \cite{cha-stri94}).\bigskip\\
 \indent In this paper we deal with martingale representations related to
 a finite set of \textit{multivariate point processes} $\{(T^i_n, X^i_n)\}_{n\geq 1},\,i=1,\cdots,d,$ defined on the same
 probability space and with explosion time taking values in  $(0,+\infty]$.~More precisely  denoting by $\mathbb{G}$ the smallest right-continuous filtration making each process of the sequence adapted (up to an initial enlargement), we deal with weak and strong representations of all $\mathbb{G}$-local martingales.~In particular we prove that all local martingales are weakly represented up to the minimum among the explosion times.~It follows that, if all explosion times are almost surely infinite, then the multidimensional semimartingale associated to the multivariate point processes satisfies the $\mathbb{G}$-WRP.~Moreover, when in addition all mark's spaces are discrete, we show that the representation also holds in strong form.~Finally we give a condition which assures that the set of local martingales driving the strong representation is a (possibly infinite) $\mathbb{G}$-basis.\bigskip\\
    \indent In order to frame our paper in the literature we observe that, when all explosion times are almost surely infinite,
our results can also be interpreted as stability of the  WRP  of a fixed component of the vector under progressive enlargement by the remaining components (see Remark \ref{stability}).~We recall that, since the classical Jeulin's book, \textit{initial enlargement} of a general filtration  $\mathbb{F}$ refers to the
    expansion of any $\sigma$-algebra of the filtration by the $\sigma$-algebra generated by a random variable, whereas \textit{progressive enlargement}
    refers in adding time by time to each element $\mathcal{F}_t$ of $\mathbb{F}$ the corresponding element of another filtration, usually called $\mathbb{H}$ (see \cite{Jeulin}).\\Several authors deal with the
    topic of propagation of martingale representations from a fixed filtration $\mathbb{F}$ to an enlargement $\mathbb{G}$.~We list here some
    papers, without claiming to be exhaustive.~As far as  the case of initial enlargement is concerned we recall the contributions
    of Amendinger \cite{ame}, Amendinger, Becherer and Schweizer \cite{ame-be-schw03}, Grorud and Pontier \cite{gro_po99} and the more recent
    paper of Fontana \cite{fontana18}.~In the framework of progressive enlargement, motivated by the applications in credit risk modeling,
    most of the authors deal with the enlargement by the natural filtration of the \textit{occurrence process of a random time} and among the others
    we mention the seminal paper of Kusuoka \cite{kusuoka99}, the contributions of Jeanblanc and Le Cam  \cite{jean-lecam-09}, Jeanblanc and
    Song \cite{jean-lecam-09}, Di Tella \cite{ditella-20}, Di Tella and Engelbert \cite{ditella-enge-pre-20}.~Some authors discuss the case when $\mathbb{H}$ is the natural
    filtration of a semimartingale (see e.g.~the papers of Xue \cite{xue}, Kchia and Protter
    \cite{kchia-prott-15} and more recently Di Tella \cite{ditella-pre-20} and Di Tella and Jeanblanc \cite{ditella-jean-21}).~In great generality in \cite{caltor15}, \cite{caltor18}, \cite{caltor-accessible-17} we studied the case when $\mathbb{H}$ is
    the reference filtration, but not necessarily the natural one, of a semimartingale.\\Here, focusing on the enlargement by multivariate point processes, we extend
    the results obtained in \cite{ditella-jean-21}  about martingale representations of the initially enlarged natural filtration of a point process progressively
    enlarged by the natural filtration of a different point process.\bigskip\\
\indent The paper is organized as follows.~In Section 2 we introduce the notations and some basic notions related to a multivariate point process.~In Section 3 we present the WRP for multivariate point processes proved by Jacod in \cite{jacod74/75}.~In Section 4 we extend Jacod's result to a finite vector of multivariate point processes.~In Section 5 we study the particular case when all explosion times are almost surely infinite and the mark's spaces are discrete.~We devote Remark \ref{confronto1 DT Jeanbl} and Remark \ref{confronto2 DT Jeanbl} to the comparison with \cite{ditella-jean-21}.
\section{Notations and basic results}
Let $(\Omega,\mathcal{F}, \mathbb{F}, P)$ be a filtered probability space with the filtration $\mathbb{F}=(\mathcal{F}_t)_{t\geq 0}$ under standard conditions.\\
\indent For $p\geq 1$, we denote by $H^p(\mathbb{F})$ the space of $p$-integrable $\mathbb{F}$-martingales, by $H^p_{0}(\mathbb{F})$ the subset of its centered elements and by $H^p_{loc}(\mathbb{F})$  and $H^p_{loc,0}(\mathbb{F})$ their localizations (see \cite{jacod} pages 26-27).~We recall that $H^1_{loc}(\mathbb{F})$ coincides with the set of all $\mathbb{F}$-local martingales (see Lemma 2.38 in \cite{jacod}).~For any two elements $M$ and $N$ of $H^1_{loc}(\mathbb{F})$ we indicate with $[M,N]$ their \textit{quadratic covariation}, that is the process defined at time $t$ by
$$[M,N]_t=\langle M^c,N^c\rangle_t+\sum_{s\leq t} \Delta M_s \Delta N_s,$$
where $M^c$ and $N^c$ are the continuous martingale parts of $M$ and $N$ respectively and $\langle M^c,N^c\rangle$ is their sharp bracket (see Definition VII.39 and Theorem VII.25 c) in \cite{del-me-b}).~$M$ and $N$ in $H^1_{loc,0}(\mathbb{F})$ are \textit{orthogonal} if and only if $MN\in H^1_{loc,0}(\mathbb{F})$.~Since $[M,N]$ is the optional compensator of $MN$, an equivalent definition of orthogonality is $[M,N]\in H^1_{loc,0}(\mathbb{F})$.\\
\indent For any measurable process $X$ we denote by $^pX$ its \textit{$\mathbb{F}$-predictable projection}, that is the predictable process uniquely defined, up to an evanescent set, by the rule
\begin{equation*}
^pX_T=E[X_T|\mathcal{F}_{T-}] \text{ on } \{T<\infty\}  \text{ for all $\mathbb{F}$-predictable times } T
\end{equation*}
(see Chapter I, Theorem 2.28 in \cite{ja-sh03}).\\
\indent For any  process $A$ with locally integrable variation we denote by $A^p$ its \textit{ $\mathbb{F}$-dual predictable projection}, that is the predictable process with locally integrable variation,  unique up to an evanescent set, such that
\begin{equation*}
A-A^p  \text{ is a $\mathbb{F}$-local martingale }
\end{equation*}
(see Chapter I, Theorem 3.18 in \cite{ja-sh03}).\\
\indent Following \cite{jacod74/75} we introduce the notion of \textit{multivariate point process} (\textit{m.p.p.}) and  its characterization by a \textit{random measure}.~We stress that this approach generalizes the more usual notion of multivariate point process, which does not include trajectories with finite explosion time (see \cite{jacod}, \cite{ja-sh03}, \cite{he-wang-yan92}).~Let $E$ be a \textit{Lusin space}, $\mathcal{E}$ its Borel $\sigma$-algebra and $\Delta$ an extra point.~Let us introduce the following notations
$$
  E_\Delta:=E\cup\{\Delta\}  \;\;\;  \tilde{E}:= (0,+\infty)\times E  \;\;\;  \tilde{E}_\Delta:= \tilde{E}\cup\{(+\infty,\Delta)\}
$$
and let us denote by $\mathcal{E}_\Delta$,
$\tilde{\mathcal{E}}$ and $\tilde{\mathcal{E}}_\Delta$ the Borel
$\sigma$-algebra of $E_\Delta$, $\tilde{E}$ and
$\tilde{E}_\Delta$, respectively.~Finally let us set
$$
\tilde{\Omega}:= \Omega\times(0,+\infty)\times E  \;\;\;  \tilde{\mathcal{P}}(\mathbb{F}):=\mathcal{P}(\mathbb{F})\otimes \mathcal{E},
$$
where $\mathcal{P}(\mathbb{F})$ is the predictable $\sigma$-algebra on $\Omega \times (0,+\infty)$.
\begin{definition}\label{def-mpp}
A multivariate point process on $(\Omega,\mathcal{F}, \mathbb{F}, P)$ is a sequence of random variables $\{(T_n,X_n)\}_{n\geq 1}$ with values in $(\tilde{E}_\Delta,\tilde{\mathcal{E}}_\Delta)$ such that
\begin{itemize}
 \item[(i)] for each $n$, $T_n$ is a $\mathbb{F}$-stopping time and $T_n\le T_{n+1}$;
 \item[(ii)] for each $n$, $X_n$ is $\mathcal{F}_{T_n}$-measurable;
 \item[(iii)] if $T_n<\infty$, then $T_n< T_{n+1}$.
\end{itemize}
\end{definition}
\begin{remark}
The value $(+\infty,\Delta)$ must be considered in order to identify any pure jump process with a m.p.p.~according to previous definition.~As an example, the  occurrence process $\mathds{1}_{[[\tau,+\infty[[}$ of a finite random time $\tau$ is a m.p.p.~with  $(T_1,X_1)=(\tau, 1)$, $(T_n,X_n)=(+\infty,\Delta)$ for any $n\geq 2$.
\end{remark}
\noindent The \textit{explosion time} of a m.p.p.~$\{(T_n,X_n)\}_{n\geq 1}$ is the random variable taking values in $(0,+\infty]$ defined by $$T_\infty:=\lim_{n}T_n.$$
\indent Any m.p.p.~$\{(T_n,X_n)\}_{n\geq 1}$  is completely characterized by  the discrete random measure $\mu$ on $(\tilde{E},\tilde{\mathcal{E}})$
defined by
\begin{equation}\label{eq-mpp-measure}
\mu(\omega; dt, dx):=\sum_{n\ge 1} \mathds{1}_{\{T_n< \infty\}}(\omega)\delta_{\left(T_n(\omega), X_n(\omega)\right)}(dt, dx),
\end{equation}
where $\delta_x$ denotes the Dirac measure concentrated on $x$.~Here \textit{random measure} means positive transition measure from $(\Omega, \mathcal{F})$ to  $(\tilde{E},\tilde{\mathcal{E}})$ (see Definition ($0.10$) in \cite{jacod}).\\
\indent Set
\begin{equation}\label{eq-G}
    \mathcal{X}_t:=\sigma\Big(\mu( (0,s]\times B): s\le t, B\in \mathcal{E}\Big).
\end{equation}
\noindent Clearly $\mathcal{X}_t\subset\mathcal{F}_t$ and the family $\mathbb{X}=(\mathcal{X}_t)_{t\ge 0}$ is a right continuous filtration, where last property depends upon the fact that the sequence $\{T_n\}_{n\geq1}$ is strictly increasing (see \cite{jacod74/75} page 238).~The filtration $\mathbb{X}$ is called the \textit{natural filtration} of $\{\left(T_n,X_n\right)\}_{n\geq1}$.\\
 \indent Let $\mu$ be the random measure defined by (\ref{eq-mpp-measure}).~Then Theorem 2.1, Proposition 2.3 and the characterization (2.6) in \cite{jacod74/75} prove existence and uniqueness of a positive random measure $\nu(\omega; dt, dx)$ on $(\tilde{E},\tilde{\mathcal{E}})$ satisfying
 \begin{align*}
    \nu(\{t\}\times E)\le 1\;\;\;\;\;\; \nu([T_\infty,\infty)\times E)=0
 \end{align*}
 such that
 \begin{itemize}
   \item[(i)] $\left(\nu(\omega; (0,t]\times B)\right)_{t\ge 0}$  is predictable for each $B\in \mathcal{E}$;
   \item[(ii)] $\left(\mu(\omega; (0,t\wedge T_n]\times B)-\nu(\omega; (0,t\wedge T_n]\times B)\right)_{t\ge 0}$ is a uniformly integrable martingale null at time zero for each $n\ge 1$ and $B\in \mathcal{E}$.
 \end{itemize}
 The measure $\nu$ is called \textit{$\mathbb{F}$-predictable compensator} or \textit{$\mathbb{F}$-dual predictable projection} of $\mu$ (see also Chapter II Theorem 1.8 page 67 in \cite{ja-sh03}).
\section{Jacod's WRP of any multivariate point process}
In this section we give a suitable version
 of a classical representation result for local martingales which are functionals of a m.p.p.~The result is due to Jacod (see \cite{jacod74/75}) and in case of infinite explosion time it coincides with the WRP of the semimartingale naturally associated to the process (see Chapter XII Section 2 Definition 13.13 page 368 in \cite{he-wang-yan92} and  %or equivalently of the multivariate point process itself
 Theorem  \ref{thm-rappr-vector2} below).
\begin{thm}\label{thm-rappr-jacod}(Theorem 5.4 in \cite{jacod74/75})\\
On a filtered probability space $(\Omega,\mathcal{F}, \mathbb{F}, P)$ let $\{(T_n,X_n)\}_{n\geq 1}$ be a m.p.p.~with values in  $(\tilde{E}_\Delta,\tilde{\mathcal{E}}_\Delta)$ and let $\mu$ be the associated random measure on $(\tilde{E},\tilde{\mathcal{E}})$ and $T_\infty$ its explosion time.
Assume that for any $t\ge 0$
\begin{equation}\label{eq-A1}
    \mathcal{F}_t=\mathcal{F}_0\vee \mathcal{X}_t
\end{equation}
 with $\mathcal{F}_0\subset\mathcal{F}$.\\
Consider a right continuous adapted process $Z=(Z_t)_{t\geq 0}$.~Then the following statements are equivalent
\begin{itemize}
  \item [(i)] there exists a sequence $\{S_n\}_{n\geq 1}$ of
  stopping times increasing a.s.~towards $T_\infty$ such that, for
  any $n\geq 1$, $Z_{t\wedge S_n}$ is a uniformly integrable
  martingale;
  \item [(ii)] there exists a finite $\tilde{\mathcal{P}}(\mathbb{F})$-measurable function $W$ on $\tilde{\Omega}$ such that on $\{t<T_\infty\}$
  \begin{equation}\label{eq-cond-rappr}
    \int_0^t\int_E|W(s,x)|\nu(ds,dx)<\infty \;\;a.s.
  \end{equation}
  and
    \begin{equation}\label{eq-jacod-rapp}
    Z_t=Z_0+\int_0^t\int_EW(s,x)(\mu(ds,dx)-\nu(ds,dx))  \;\;a.s.
  \end{equation}
\end{itemize}
\end{thm}
\begin{remark}\label{rem:explosion}
Following \cite{jacod74/75} we can express
point $(i)$ above saying that $Z$ is a right continuous $\mathbb{F}$-local martingale on $[0,T_\infty)$ and writing $Z\in \mathcal{M}^{T_\infty}_{loc}(\mathbb{F})$.~It is worthwhile to observe that
$$H^1_{loc}(\mathbb{F})\subset\mathcal{M}^{T_\infty}_{loc}(\mathbb{F})$$
  and when $P(T_\infty<+\infty)=0$ the sets $H^1_{loc}(\mathbb{F})$ and $\mathcal{M}^{T_\infty}_{loc}(\mathbb{F})$ coincide.~In fact, if $Z$ is a right continuous $\mathbb{F}$-local martingale with reducing sequence $(\tau_n)_{n\geq 1}$, then the sequence of stopping times $(\tau_n\wedge T_n)_{n\geq 1}$ plays the role of the sequence $(S_n)_{n\geq 1}$ of condition (i).\\
Then in particular previous theorem states  that all $\mathbb{F}$-local martingales can be weakly represented up to $T_\infty$.
\end{remark}
\noindent In the light of the above remark it comes natural to express Jacod's representation theorem saying that\\
\noindent \textit{any m.p.p.~satisfies the WRP  up to $T_\infty$ (just WRP if $T_\infty=+\infty$ $a.s.$) with respect to its initially enlarged natural filtration.}\bigskip\\
\indent As an immediate consequence recalling the definition of WRP of semimartingales  we get the following result (see Definition 13.13 Section 2 Chapter XII page 368 in \cite{he-wang-yan92}).
\begin{cor}\label{weak prp of m.p.p.}
When $P(T_\infty<+\infty)=0$ the semimartingale $(X_t)_{t\geq 0}$ defined by
\begin{equation}\label{eq-mpp-processes}
    X_t:=\sum_{n\ge 1}X_n\mathds{1}_{\{T_n\le t\}}
\end{equation}
 satisfies the $\mathbb{F}$-WRP.
 \end{cor}
 \begin{proof}
Let us observe that, according to Chapter XI Definition 11.1 page 293 in \cite{he-wang-yan92}, we need to extend the random measure $\mu$ defined by (\ref{eq-mpp-measure}) on $(0,+\infty)\times E$, to the measurable space given by $[0,+\infty)\times E$ endowed with its Borel $\sigma$-algebra.~Such an extension coincides with the jump measure of $(X_t)_{t\geq 0}$ and the analogous extension of $\nu$ coincides with the L\'evy system of $(X_t)_{t\geq 0}$ (see Chapter XI Theorem 11.15 page 300 in \cite{he-wang-yan92}).~Then the thesis immediately follows by Theorem \ref{thm-rappr-jacod}.
 \end{proof}
\section{Propagation of weak martingale representation under enlargement by multivariate point processes}\label{sec 2}
In this section we
point out how from the  WRP of any
m.p.p.~easily follows the WRP of a finite vector of m.p.p.s.~By the WRP of a finite vector of m.p.p.s.~we mean that any local martingale with respect to the union of the initially enlarged vector's natural filtration is represented as integral of a predictable function with respect to a unique compensated random measure.~More precisely the representation holds up to the minimum among the explosion times of the components.\\
Our result
   implies the propagation of martingale representation from the initially enlarged natural filtration  of the first component of the considered vector  to its progressive enlargement by the union of the initially enlarged natural filtrations  of the remaining  components.\\

\indent Let $(\Omega,\mathcal{F}, P)$ be a probability space supporting  $d$ filtrations  under standard conditions, $\mathbb{F}^i,\, i=1,\ldots , d$.~For any $i$ let $\{(T^i_n, X^i_n)\}_{n\geq 1}$  be a m.p.p.~on $(\Omega,\mathcal{F}, \mathbb{F}^i, P)$ with values in $(\tilde{E}^i_\Delta,\tilde{\mathcal{E}}^i_\Delta)$ with   $E^{i}$ a Lusin space.~Assume that for any $t\ge 0$
\begin{equation*}
    \mathcal{F}^i_t=\mathcal{F}^i_0\vee \mathcal{X}^i_t
\end{equation*}
 where $\mathcal{F}^i_0\subset\mathcal{F}$ and $\mathbb{X}^i=(\mathcal{X}^i_t)_{t\ge 0}$ is the natural filtration of $\{(T^i_n, X^i_n)\}_{n\geq 1}$.~Set  $$E^i_0:=E^i\cup\{0\},$$
 where, with a little abuse of language, the same symbol $0$  denotes the neutral element of the compact space containing $E^i$ for any $i$.~Finally let $\mathbb{G}=(\mathcal{G}_t)_{t\ge 0}$ be the filtration defined by
 \begin{equation}\label{eq newG}
\mathcal{G}_t:=\cap_{s>t}\,\vee_{i=1}^{d}\mathcal{F}^i_s.
 \end{equation}
  \bigskip
\noindent The natural candidate for the $\mathbb{G}$-WRP is the \textit{merging process} of the $d$ assigned m.p.p.s up to the minimum explosion time, that is a new m.p.p.~$\{(T_n, V_n)\}_{n\geq 1}$ on $(\Omega,\mathcal{F}, \mathbb{G}, P)$ taking values in $(\tilde{E}_\Delta,\tilde{\mathcal{E}}_\Delta)$, with
    \begin{equation}\label{eq-vector-space}
    E:=E^1_0\times E^2_0 \ldots \times E^d_0,
    \end{equation}
 and having explosion time $T_\infty$ such that
 \begin{equation}\label{explosion time}T_\infty=\min \left(T^1_\infty,\ldots,T_\infty^d\right).\end{equation}
Let us be more precise.\\First of all we construct  recursively the sequence of random times $\{T_n\}_{n\geq 1}$.\\

\noindent $T_1(\omega):=\begin{cases}\inf\{T^i_1(\omega):\; T^i_1(\omega)<+\infty,\ \ i=1,...,d\} & \text{ if } \{\ldots\}\neq\emptyset\\
  +\infty & \text{  otherwise }
     \end{cases}$\\
$T_2(\omega):=\begin{cases}\inf\{T^i_k(\omega):\; T_1(\omega)<T^i_k(\omega)<+\infty ,\ \ i=1,...,d,\;k\ge 1\} & \text{ if } \{\ldots\}\neq\emptyset\\
  +\infty & \text{  otherwise }
     \end{cases}$\\
\vdots\\
$T_n(\omega):=\begin{cases}\inf\{T^i_k(\omega):\; T_{n-1}(\omega)<T^i_k(\omega)<+\infty ,\ \ i=1,...,d,\;k\ge 1\} & \text{ if } \{\ldots\}\neq\emptyset\\
  +\infty & \text{  otherwise }
     \end{cases}$\\
 \vdots\\
  It is easy to check that the sequence $\{T_n\}_{n\geq 1}$ satisfies properties (i) and (iii) of  Definition \ref{def-mpp} with respect to  $\mathbb{G}$.\\
Next we construct the sequence of random variables $\{V_n\}_{n\geq 1}$.
 \begin{equation*}
    V_n:=\begin{cases} (V^1_n,\cdots , V^d_n) & \text{ if } T_n<+\infty  \\
   \Delta & \text{  otherwise }
     \end{cases}
     \end{equation*}
 where, when $T_n<+\infty$, the $i$-component $V^i_n$ is the random variable defined by
 \begin{equation*}
    V^i_n(\omega):=\begin{cases}
    X^i_{k}(\omega) & \text{  if there exists  }  k\geq 1  \text{ such that  }   {T_k^i}(\omega)=T_n(\omega)\\
  0 & \text{  otherwise. }
    \end{cases}
 \end{equation*}
  By definition $\{V_n\}_{n\geq 1}$ takes values in $E \cup \{\Delta\}$ and
  \begin{equation*}
    V^i_n=\sum_{k\ge 1}X^i_{k}\mathds{1}_{\{{T_k^i}=T_n\}}\mathds{1}_{\{T_n<+\infty\}}+\Delta \, \mathds{1}_{\{T_n=+\infty\}}.
 \end{equation*}
 \noindent Since the random time  $T_n$ and those in the sequence $\{T_k^i\}_{k\ge 1}$ are $\mathbb{G}$-stopping times, then by Chapter I Property 1.17 page 4 in \cite{ja-sh03} it follows
   $$
   \mathcal{G}_{T^i_k}\cap {\{{T_k^i}=T_n\}}=\mathcal{G}_{T_n}\cap {\{{T_k^i}=T_n\}}.
   $$
By definition $X^i_{k}$ is  $\mathcal{G}_{T^i_k}$-measurable so that previous equality implies that the random variable $X^i_{k}\mathds{1}_{\{{T_k^i}=T_n\}}$ is  $\mathcal{G}_{T_n}$-measurable and therefore  for each $n\geq 1$, $V_n$ satisfies (ii) of Definition (\ref{def-mpp})  with respect to $\mathbb{G}$.\\
Finally by construction $T_\infty\leq T^i_{\infty}$ for all $i$, so that (\ref{explosion time}) holds.~Therefore the sequence  $\{(T_n, V_n)\}_{n\geq 1}$ is the m.p.p.~we were looking for.
 \begin{remark}\label{set no zero}
Note that the random variables $V_n,\,n=1,2,\ldots,$ do not assume the value $\{0,\ldots,0\}$, which could therefore be excluded from the mark's space.~Moreover,  since in general the m.p.p.s.~do not share all random times, the new mark's space could be strictly contained in $E\setminus \{0,\ldots,0\}$.~Indeed $E\setminus \{0,\ldots,0\}$
can be represented as
 \begin{equation*}\label{set E} \cup_{k=1}^{d}\,\cup_{i_1<\ldots <i_k} E^{i_1,\ldots,i_k},\end{equation*}
 where for any $k\in \{1,\ldots,d\}$ and for any choice of $(i_1,\ldots ,i_k)\in\{1,\ldots,d\}^k$ such that $i_1<\ldots <i_k$%\,i_l\in\{1,\ldots,d\}$
 $$E^{i_1,\ldots,i_k}:=\{x=(x_1,\ldots,x_d)\in E^1_0\times \ldots \times E^d_0 :\, x_{i}=0,\;\text{iff}\;\; i\notin (i_1,\ldots,i_k)\}$$
  is the Lusin space, which corresponds to the simultaneous jumps of the processes indexed by the set $\{i_1,\ldots,i_k\}$ while the others stay fixed.~In particular, when $k=1$, for any choice of $i_1$, $E^{i_1}$ tells that the process indexed by $i_1$ is jumping alone and at the opposite, when $k=d$, $E^{i_1,\ldots,i_d}=E^{1,\ldots,d}=E^1\times\ldots\times E^d$ tells that all processes are jumping together.
 \end{remark}
 \bigskip
 Let us now state the main theorem of this section.
 \begin{thm}\label{thm-rappr-vector1}
In the framework of this section let $\mathbb{X}=(\mathcal{X}_t)_{t\geq 0}$ be the natural filtration of $\{(T_n, V_n)\}_{n\geq 1}$.~Set
 \begin{equation}\label{F_0}\mathcal{F}_0:=\vee_{i=1}^{d}\mathcal{F}^i_0.\end{equation}
Let $\mathbb{F}=(\mathcal{F}_t)_{t\geq 0}$ be the right-continuous filtration defined by
   \begin{equation}\label{filtration F2}
 \mathcal{F}_t:=\mathcal{F}_0\vee \mathcal{X}_t.
 \end{equation}
Then  $\{(T_n, V_n)\}_{n\geq 1}$ satisfies the $\mathbb{F}$-WRP up to $T_\infty.$
 \end{thm}
 \begin{proof}
 The right-continuity of the filtration $\mathbb{F}$ derives immediately by the same property of $\mathbb{X}$ and this in turn easily follows by Theorem 23 (iii) in Appendix 2 of \cite{brem-81} where the following characterization is proved
\begin{equation*}
\mathcal{X}_t=\sigma\Big(\mathds{1}_{\{V_n\in B\}}\mathds{1}_{\{T_n\le s\}}:  n\ge 1, s\le t, B\in \mathcal{E}\Big).
\end{equation*}
  Then  $\{(T_n, V_n)\}_{n\geq 1}$ is a m.p.p.~on $(\Omega,\mathcal{F}, \mathbb{F}, P)$ and
 Theorem \ref{thm-rappr-jacod} applies to it, since (\ref{filtration F2}) holds and $\mathcal{F}_0\subset \mathcal{F}$.\end{proof}
\begin{lemma}\label{lemma-no explosion}
 Assume  $P(T^i_\infty<+\infty)=0$, for each $i=1,...,d$.~Then
 \begin{equation*}\mathbb{F}=\mathbb{G},\end{equation*}
 where $\mathbb{G}$ is the filtration defined by (\ref{eq newG}).
\end{lemma}
\begin{proof}
  Recalling  (\ref{filtration F2}), (\ref{F_0}) and (\ref{eq newG}) it sufficies to show that
  $$\mathcal{X}_t=\cap_{s>t}\,\vee_{i=1}^{d}\mathcal{X}^i_s.$$
 To this end we start by observing that the assumption on the explosion times implies $P(T_\infty<+\infty)=0$.~Since the process $\{(T_n, V_n)\}_{n\geq 1}$ is a function of the set of processes $\{(T^i_n, X^i_n)\}_{n\geq 1}, i=1,\ldots,d$, then the inclusion $\mathcal{X}_t
 \subset \bigcap_{s>t}\,\vee_{i=1}^{d}\mathcal{X}^i_s$ is immediate.~In order to prove the reverse inclusion, is it enough to prove that $\mathcal{X}^i_t
 \subset\mathcal{X}_t$ for any $t\in \mathbb{R}^+$ and for each $i=1,\ldots,d$.~We take w.l.o.g.~$i=1$ and we prove that for any $t>0$
 \begin{equation}\label{inclusion}\mathcal{X}^1_t\subset\mathcal{X}_t.\end{equation}
As in (\ref{eq-G}) it holds
\begin{equation*}
    \mathcal{X}^1_t:=\sigma\Big(\mu^{1}( (0,s]\times B^1): s\le t, B^1\in \mathcal{E}^1\Big).
\end{equation*}
 where $\mu^{1}$ is the random measure associated to $\{(T^1_n, X^1_n)\}_{n\geq 1}$ by  (\ref{eq-mpp-measure}).~Therefore, if we prove that for each $s\leq t$ and $B^1\in \mathcal{E}^1$
\begin{equation}\label{reverse}\mu^{1}( (0,s]\times B^1)=\mu\left( (0,s]\times (B^1 \times E^2_0 \times \ldots E^d_0)\right), \  \ \text{a.s.},  \end{equation}
then we get (\ref{inclusion}).~In fact
 \begin{align*}
   \mu\left( (0,s] \times (B^1 \times E^2_0 \times \ldots E^d_0)\right)=\sum_{n\ge 1}\mathds{1}_{\{T_n\le t\}}\mathds{1}_{\{V^1_n\in B^1\}}
 \end{align*}
or equivalently since $0\notin B^1$
 \begin{align*}
   \mu\left ((0,s]\times (B^1 \times E^2_0 \times \ldots E^d_0)\right)=\sum_{n\ge 1}\mathds{1}_{\{T_n\le t\}}\mathds{1}_{\{V^1_n\in B^1\}}\mathds{1}_{\{V^1_n\neq 0\}}.
 \end{align*}
Now, by construction
$\mathds{1}_{\{V^1_n\neq 0\}}=\sum_{k\ge 1}\mathds{1}_{\{T_n=T^1_k\}}\mathds{1}_{\{V^1_n=X^1_k\}}$
so that
 \begin{align*}
   &\mu\left((0,s]\times (B^1 \times E^2_0 \times \ldots E^d_0)\right)=\cr
   &\sum_{n\ge 1}\mathds{1}_{\{T_n\le t\}}\mathds{1}_{\{V^1_n\in B^1\}}\sum_{k\ge 1}\mathds{1}_{\{T_n=T^1_k\}}\mathds{1}_{\{V^1_n=X^1_k\}}=\cr
   &\sum_{n\ge 1} \sum_{k\ge 1}\mathds{1}_{\{T^1_k\le t\}}\mathds{1}_{\{X^1_k\in B^1\}}\mathds{1}_{\{T_n=T^1_k\}}\mathds{1}_{\{V^1_n=X^1_k\}}=\cr
   &\sum_{k\ge 1}\mathds{1}_{\{T^1_k\le t\}}\mathds{1}_{\{X^1_k\in B^1\}}\sum_{n\ge 1} \mathds{1}_{\{T_n=T^1_k\}}\mathds{1}_{\{V^1_n=X^1_k\}}.
    \end{align*}
    Since $P(T_\infty<+\infty)=0$ then
   $$\sum_{n\ge 1}\mathds{1}_{\{T_n=T^1_k\}}\mathds{1}_{\{V^1_n=X^1_k\}}=1, \  \ \text{a.s.}$$
   and therefore
    \begin{align*}
     \mu\left((0,s]\times (B^1 \times E^2_0 \times \ldots E^d_0)\right)=
    \sum_{k\ge 1}\mathds{1}_{\{T^1_k\le t\}}\mathds{1}_{\{X^1_k\in B^1\}}, \  \ \text{a.s.}
 \end{align*}
 that is equality (\ref{reverse}).
\end{proof}
Theorem \ref{thm-rappr-vector1} and Lemma \ref{lemma-no explosion} easily yield next result.
\begin{thm}\label{thm-stability}
Assume  $P(T^i_\infty<+\infty)=0$, for each $i=1,...,d$.~Then
 the m.p.p.~$\{(T_n, V_n)\}_{n\geq 1}$ satisfies the $\mathbb{G}$-WRP.
\end{thm}
 \begin{remark}\label{stability}
   It is clear that under the assumption of the theorem the right-continuous version of $\mathbb{F}^2\vee\ldots\vee\mathbb{F}^d$  coincides with the right-continuous version of the natural filtration of the merging of $\{(T^2_n, X^2_n)\}_{n\geq 1}, \ldots, \{(T^d_n, X^d_n)\}_{n\geq 1}$ initially enlarged by $\mathcal{F}^2_0\vee.....\vee\mathcal{F}^d_0$.~Then the above theorem can be interpreted as stability of the WRP of $\{(T^1_n, X^1_n)\}_{n\geq 1}$ under progressive enlargement by
 $\{(T^2_n, X^2_n)\}_{n\geq 1}, \ldots, \{(T^d_n, X^d_n)\}_{n\geq 1}$.
 \end{remark}
\begin{remark}
In the frame of the theory of stable subspaces, when $T_\infty=+\infty$,
following Definition  4.45 and Theorem 4.46 in \cite{jacod},
 Theorem \ref{thm-stability} implies that
\begin{equation*}
H^1_0(\mathbb{G})\equiv K^{1,1}(\mu),
\end{equation*}
where $K^{1,1}(\mu)$ is the stable subspace of $H^1_0(\mathbb{G})$ defined by
$$K^{1,1}(\mu)=\left\{\int_0^\cdot\int_EW(s,x)(\mu(ds,dx)-\nu(ds,dx)), \ W \in G^1(\mu)\right\},$$
 with  $G^1(\mu)$ as in   Definition 3.62 at page 98 in \cite{jacod}.\\
 In fact, trivially  $K^{1,1}(\mu)\subset H^1_0(\mathbb{G})$.~To prove the reverse inclusion, observe that if $Z\in H^1_0(\mathbb{G})$, then $Z\in H^1_{loc,0}(\mathbb{G})$ so by Theorem \ref{thm-rappr-vector1} there exists $W$  finite  and  $\tilde{\mathcal{P}}$-measurable satisfying (\ref{eq-cond-rappr}) and (\ref{eq-jacod-rapp}).~Therefore, to prove that $Z\in K^{1,1}(\mu)$  it is enough  to prove that  $W\in G^1(\mu)$, that is
 \begin{equation*}
 E\left[\left(\sum_{s} \left(\int_EW(s,x)(\mu(\{s\},dx)-\nu(\{s\},dx))\right)^2\right)^\frac{1}{2}\right]<+\infty.
 \end{equation*}
 This fact derives by noting that, since $Z$ is locally of finite variation, then
 \begin{equation*}
 [Z,Z]^\frac{1}{2}_\infty=\left(\sum_{s} \Delta Z^2_s\right)^\frac{1}{2}=\left(\sum_{s} \left(\int_EW(s,x)(\mu(\{s\},dx)-\nu(\{s\},dx))\right)^2\right)^\frac{1}{2}
 \end{equation*}
 and, since $Z\in H^1_0(\mathbb{G})$, then $E[Z^*_\infty] <+\infty$  and moreover
  \begin{equation*}
 E[Z^*_\infty]\le a E[[Z,Z]^\frac{1}{2}_\infty]\le b E[Z^*_\infty]
 \end{equation*}
 (see 2.1 and 2.34 in \cite{jacod}).
\end{remark}
 \begin{thm}\label{thm-rappr-vector2}
 In the setting of this section assume $P(T^i_\infty<+\infty)=0$, for each $i=1,...,d$.~Let $(X^i_t)_{t\geq 0}$ be the semimartingale associated to $\{(T^i_n, X^i_n)\}_{n\geq 1}$ by the analogous of (\ref{eq-mpp-processes}).~Then the semimartingale $(X_t)_{t\geq 0}$ defined by
 $$X_t:=(X^1_t,\ldots,X^d_t)$$
 satisfies the $\mathbb{G}$-WRP.
\end{thm}
\begin{proof}
Since $P(T_\infty<+\infty)=0$,
\begin{equation*}
    V_t:=\sum_{n\ge 1}V_n\mathds{1}_{T_n\le t}
\end{equation*}
is well posed $a.s.~$for all $t\geq 0$.~Corollary \ref{weak prp of m.p.p.} applies to the semimartingale $(V_t)_{t\geq 0}$
 and says that it satisfies the $\mathbb{G}$-WRP.\\
 Let us prove that $(V_t)_{t\geq 0}$ and $(X_t)_{t\geq 0}$ are indistinguishable processes and more precisely, if we set $$N:=\{T_\infty<+\infty\},$$ then when $\omega\in N^c$ for all $t\geq 0$
\begin{equation*}
V_t(\omega)=X_t(\omega).
\end{equation*}
In fact, fixed $i=1,\ldots,d,$ on $N^c$
 \begin{align*} V^i_t&=\sum_{n\ge 1}V^i_n\mathds{1}_{\{T_n\le t\}}
 \end{align*}
 and, considering that the general term of the sum is different from 0 only if there exists $k$ such that ${T^i_k=T_n}$ and $V^i_n=X^i_k$, on $N^c$
 \begin{align*} V^i_t&=\sum_{n\ge 1}V^i_n\mathds{1}_{\{T_n\le t\}}\sum_{k\ge 1}\mathds{1}_{\{T^i_k=T_n\}}\cr&=
 \sum_{n\ge 1}\sum_{k\ge 1}X^i_k\mathds{1}_{\{T^i_k\le t\}}\mathds{1}_{\{T^i_k=T_n\}}.
 \end{align*}
Due to the fact that on $N^c$  the sum is finite, the order of summation can be changed providing
 \begin{equation*}
 V^i_t=\sum_{k\ge 1}X^i_k\mathds{1}_{\{T^i_k\le t\}}\sum_{n\geq 1}\mathds{1}_{\{T^i_k=T_n\}}.
   \end{equation*}
   By construction on $N^c$  $$\sum_{n\geq 1}\mathds{1}_{\{T^i_{k}=T_n\}}=1$$
    and therefore
  \begin{align*}
   V^i_t&=\sum_{k\ge 1}X^i_k\mathds{1}_{\{T^i_k\le t\}}= X^i_t.
  \end{align*}
We conclude the proof observing that the indistinguishability of the processes $(V_t)_{t\geq 0}$ and $(X_t)_{t\geq 0}$ implies that their jump measures coincide $a.s.$
\end{proof}
\begin{remark}\label{confronto1 DT Jeanbl}
Point (i) of Theorem 3.5 in \cite{ditella-jean-21} is a particular case of previous theorem.~More precisely in \cite{ditella-jean-21} the authors deal with a pair of non explosive point processes, a particular case of our setting, that is
when $d=2$, $E^1=E^2=\{1\}$ and hence $E\setminus \{0,0\}=\{(0,1), (1,0), (1,1)\}$.~They identify the two point processes with the jump measure of two different semimartingales, $X$ and $H$.~They call $\mathbb{F}$ the initially enlarged natural filtration of $X$ and $\mathbb{H}$ the natural filtration of $H$ and finally they indicate with $\mathbb{G}$ the right-continuous version of  $\mathbb{F}\vee\mathbb{H}$.~They show that the semimartingale $(X,H)$ has the $\mathbb{G}$-WRP.
\end{remark}
\section{The non explosive case with discrete marks' spaces}
In this section, working in the same setting and with the same notations of Section \ref{sec 2},
 we assume that for each $i\in\{1,...,d\}$ the explosion time $T^i_\infty$ is infinite a.s.~and the marks' space
$E^i$ is discrete, so that the same happens for $T_\infty$ and $E$ defined by (\ref{explosion time}) and (\ref{eq-vector-space}) respectively.\bigskip\\
\indent The first goal is to prove that there exists a discrete family of local martingales which strongly represents
 $\mathbb{G}$ (see (\ref{eq newG})).~The idea is that the multidimensional m.p.p.~$\{(T_n, V_n)\}_{n\geq 1}$, which as stated in Theorem \ref{stability} weakly represents any $\mathbb{G}$-local
martingale, can be \textit{splitted}
in a discrete family of point processes, one for each point of the mark's space $E$, except the null element $\{0,\ldots,0\}$.~Trivially the semimartingales
associated to these point processes are increasing processes, do not share any jump time and each of them  defines, up predictable compensation, a
local martingale.~Then under a suitable assumption we show that the above set of local martingales  is  a $\mathbb{G}$-basis.
\subsection{From WRP to PRP}
First of all let us construct the family of local martingales strongly representing $\mathbb{G}$.~We are in the framework of Theorem \ref{thm-stability}, so that $\{(T_n, V_n)\}_{n\geq 1}$ satisfies the $\mathbb{G}$-WRP.~Let $\mu$ be the random measure associated to $(T_n,V_n)$ by formula (\ref{eq-mpp-measure})
 and let $\nu$ be its $\mathbb{G}$-dual predictable projection.~Then for any $Z\in H^{1}_{loc}(\mathbb{G})$
\begin{align*}
Z_t = Z_0 + \int_0^t\int_E W(s,x)\big(\mu(ds,dx)-\nu(ds,dx)\big)
,\;\;\;a.s.
\end{align*}
where $W$ is a
function  as in point (ii) of Theorem \ref{thm-rappr-jacod}.~By construction $\mu-\nu$ is a random measure on the discrete set $E\setminus \{0,\ldots,0\}$ (see Remark \ref{set no zero}).~Without loss of generality we consider $E$ countable and we set
$$E\setminus \{0,\ldots,0\}=\{x_1,x_2,\ldots\}.$$
Then we can rewrite
\begin{align*}
Z_t =& Z_0 +
\int_0^t\int_{E\setminus \{0,\ldots,0\}}\sum_{h\geq 1}\,W(s,x)\mathds{1}_{\{x=x_h\}}\big(\mu(ds,dx)-\nu(ds,dx)\big)=\\
=& Z_0 + \sum_{h\geq 1}\,\int_0^t\,W(s,x_h)\big(\mu(ds,\{x_h\})-\nu(ds,\{x_h\})\big),\;\;\;a.s.
\end{align*}
\indent Let us set
\begin{equation}\label{M^h}
M^h_t:=\mu((0,t],\{x_h\})-\nu((0,t],\{x_h\})
\end{equation}
and
\begin{equation*}
W_t(x_h):=W(t,x_h).
\end{equation*}
We observe that the definition of $\nu$ implies that, for any fixed $h\geq 1$, the process $\left(M^{h}_t\right)_{t\geq 0}$
is a $\mathbb{G}$-local martingale null at $0$ and the $\widetilde{\mathcal{P}}(\mathbb{G})$-measurability of $W$
implies that the process  $\left(W_t(x_h)\right)_{t\geq 0}$ is
a $\mathbb{G}$-predictable process.~Then for each $t\geq 0$
\begin{align}\label{representation 1}
Z_t=Z_0 + \sum_{h\geq 1}\,\int_0^t\,W_s(x_h)dM^h_s,\;\;\;a.s.
\end{align}
Any Lebesgue-Stieltjes integral in the above formula makes sense by the integrability condition satisfied by $W$ (see point (ii) of Theorem \ref{thm-rappr-jacod}).~Moreover it can be easily proved that for any $h$ the process $\left(W_t(x_h)\right)_{t\geq 0}\in L^1_{loc}(M^h)$ that is $\left(\int _0^t\,W^2_s(x_h)d[M^h]_s\right)^\frac{1}{2}_{t\geq 0}$ is a locally integrable process (see (2.40) page 42 in \cite{jacod}).~In fact
$$\Big(\int _0^t\,W^2_s(x_h)d[M^h]_s\Big)^{1/2}=\Big(\sum_{s\leq t}(W_s(x_h)\Delta M^h_s)^2\Big)^{1/2}\leq \sum_{s\leq t}|W_s(x_h)\Delta M^h_s|$$
(see point 1 of Theorem 9.5 in \cite{he-wang-yan92}).~Therefore any integral in the sum coincides with a stochastic integral
(see Definition 2.46 and Proposition 2.48 page 45 in \cite{jacod}).\\

\indent Set  $$M:=(M^1,\ldots,M^h,\ldots)$$ with $M^h$ defined by (\ref{M^h}).~Because of the arbitrariness of $Z$ in equality (\ref{representation 1}), $\mathbb{G}$ is strongly represented by $M$ in the sense that
for any $Z\in H^{1}_{loc}(\mathbb{G})$
\begin{align*}
Z=Z_0 + \sum_{h\geq 1}\,W(x_h)\cdot M^h
\end{align*}
and we get the following result.\\

\begin{proposition}\label{prop-strong G}
$M$ enjoys the $\mathbb{G}$-PRP.
\end{proposition}

\begin{remark}\label{stable subspace 2}
Let $p\geq 1$.~Observe that $M$ is a subset of $H^p_{loc,0}(\mathbb{G})$.~Therefore under the assumptions of this section, by the same techniques used for proving point (i) of Corollary 3.6 in \cite{ditella-jean-21}, it can be shown that
 $$H^p_{0}(\mathbb{G})=\mathcal{Z}^p(M),$$
 where $\mathcal{Z}^p(M)$ is the stable subspace generated by $M$ in $H^p(\mathbb{G})$ (see Definition 4.4, Proposition 4.5 and Theorem 4.6 in \cite{jacod}).
\end{remark}
\subsection{A sufficient condition for the orthogonality  }
Consider now the family of  $\mathbb{G}$-point processes $N^{h},\,h\geq 1,$ defined by
\begin{equation*}
N^{h}_t:=\mu((0,t],\{x_h\}).
\end{equation*}
 If we denote by $N^{h,p}$ the $\mathbb{G}$-dual predictable projection of $N^{h}$, then, by its uniqueness
 and by the definition of $\nu$, we derive
\begin{equation*}
\nu((0,t],\{x_h\})=N^{h,p}_t,
\end{equation*}
so that we can write
\begin{equation}\label{eq-martingales}
M^h=N^{h}-N^{h,p}.\end{equation}
\indent In order to give a sufficient condition under which $\{M^{h}\}_{h\geq 1}$  is
 a sequence of pairwise orthogonal $\mathbb{G}$-local martingales  we need
two preliminary results.
\begin{lemma}\label{lemma-jumps}
Let $X$ be a $\mathbb{G}$-adapted pure jump process, that is such that
\begin{equation*}
X_t=\sum_{s\leq t}\,\Delta X_s
\end{equation*}
 for all $t$.~Let us assume that $X$ is of locally integrable variation and let $X^p$ be its $\mathbb{G}$-dual predictable projection.~Then, for any fixed  stopping time $S$,
$\Delta X_S=0$, a.s., implies $\Delta X^p_S=0$, a.s.
\end{lemma}
\begin{proof}
First of all note that $X$ is of finite variation
since
$$
X_t=\sum_{s\le t}\Delta X_s\mathds{1}_{(0,+\infty)}(\Delta
X_s)-\sum_{s\le t}|\Delta X_s|\mathds{1}_{(-\infty,0)}(\Delta X_s)
$$
and recall that  a finite variation process is locally integrable if and only if  there exists the corresponding dual predictable  projection  (see Definition VI-79  and Theorem VI-80 in
\cite{del-me-b}).~The proof can be easily reduced to the case when $S$ is a totally inaccessible or an accessible random time.~Now if $S$ is total inaccessible then $\Delta X^p_S=0$, a.s., since a predictable random process cannot jump on totally inaccessible times
(see Proposition 2.44 in \cite{nike06}).~If $S$ is accessible then there exists a sequence $\{S_n\}_{n\geq 1}$ of predictable jump times enveloping $S$.~Then by VI-(81.2) in \cite{del-me-b}, for each element $S_n$ of such a sequence
 $$E[\Delta X_{S_n}|\mathcal{G}_{S_n^-}]=\Delta X^p_{S_n}.$$
 This in particular implies that if $\Delta X^p_{S_n}\neq 0$ with positive probability then it has necessary to be
 $\Delta X_{S_n}\neq 0$ with positive probability.~It follows that if $\Delta X_{S}=0$ with probability one the same happens for its predictable compensator (see also Lemma 5.41 page 156 in \cite{he-wang-yan92}).
  \end{proof}
 Let us introduce for two general $\mathbb{G}$-adapted locally integrable pure jump processes $X$ and $Y$ the following condition, which we call of
\textit{mutual avoiding predictable jump times}.
\begin{ass}\label{ass:NOPJT}
$P(\Delta X_{\sigma}\neq 0)>0$ implies $\Delta Y_{\sigma}=0$, a.s.,~for any finite $\mathbb{G}$-predictable stopping time $\sigma$.
\end{ass}
\begin{remark}\label{on NOPJT} It is to stress that previous assumption is symmetric in $X$ and $Y$.~Note moreover that the condition is well-posed, since for example it holds either when one of the two processes admits $\mathbb{G}$-totally inaccessible jump times only or when $Y=\mathds{1}_{[[\tau,+\infty[[}$ with $\tau$ a finite random time which avoids the jump times of $X$.
\end{remark}
  \begin{proposition}\label{prop-orth}
Let $X$ and $Y$ be $\mathbb{G}$-adapted pure jump
processes of locally integrable variation
which verify Assumption \ref{ass:NOPJT}.~Then for any finite $\mathbb{G}$-predictable stopping time $\sigma$
$$\Delta X^p_{\sigma}\Delta Y^p_{\sigma}=0,\;\;\; \textrm{a.s.}$$
\end{proposition}
\begin{proof}
Let us suppose the existence of a $\mathbb{G}$-predictable finite stopping time $\sigma$ such that
$$P(\Delta X_{\sigma}\neq 0)>0.$$
Then Assumption \ref{ass:NOPJT} implies $\Delta Y_{\sigma}=0$, a.s., and
Lemma \ref{lemma-jumps} applies to prove that $\Delta Y^p_{\sigma}=0$, a.s., and therefore
\begin{equation}\label{intermediate}\Delta X_\sigma\Delta Y^p_\sigma= 0,  \ \ \text{a.s.}\end{equation}
Property 3.21 at page 33 in \cite{ja-sh03} assures that
$$\Delta Y^p_\sigma=\prescript{p}{}(\Delta Y)_\sigma,$$
where $\prescript{p}{}(\Delta Y)$ denotes the predictable projection of the process given by the jumps of $Y$ (see page 23 in \cite{ja-sh03}),  and then
$$\Delta X_\sigma  \prescript{p}{}(\Delta Y)_\sigma= 0,  \ \ \text{a.s.}$$
Taking the predictable projection (see Theorem 2.28 in \cite{ja-sh03}) we obtain
$$^p\left(\Delta X \prescript{p}{}(\Delta Y)\right)_\sigma=0,  \ \ \text{a.s.}$$
 and previous equality corresponds to the thesis since Remark VI-44 (e) in
\cite{del-me-b} implies
$$^p\left(\Delta X \prescript{p}{}(\Delta Y)\right)_\sigma=\prescript{p}{}(\Delta X)_\sigma  \prescript{p}{}(\Delta Y)_\sigma$$
and using again Property 3.21 in \cite{ja-sh03}
$$^p\left(\Delta X \prescript{p}{}(\Delta Y)\right)_\sigma=
\Delta X^p_\sigma\Delta Y^p_\sigma$$
(see also Theorem 2.28 at page 23 in \cite{ja-sh03}).
\end{proof}
\begin{cor}\label{rem:1}
In the setting of Proposition \ref{prop-orth} $$\Delta X_\sigma\Delta Y^p_\sigma=\Delta X^p_\sigma\Delta Y_\sigma=0,  \ \ \text{a.s.}$$
\end{cor}
\begin{proof}
In fact (\ref{intermediate}) holds and, as already recalled, Assumption \ref{ass:NOPJT} is symmetric in  $X$ and $Y$.
\end{proof}
\begin{remark}
 The condition $\Delta X_\sigma\Delta Y_\sigma= 0$, a.s., for all finite $\mathbb{G}$-predictable stopping time $\sigma$, which is weaker than Assumption \ref{ass:NOPJT}, does not prevent that $X^p$ and $Y^p$ jump together with positive probability.~Counterexample A.2.~in \cite{ditella-jean-21} provides a pair of point processes without common jumps, whose dual predictable projections share a
 jump with probability one.~It is easy to check that those processes do not satisfy our Assumption \ref{ass:NOPJT}.
\end{remark}
We are now able to state the announced result of representation for $\mathbb{G}$.
\begin{thm}\label{thm-ort-basis}
 In the setting of this section let Assumption (\ref{ass:NOPJT}) be in force for all pair $N^{h},N^{k}$ with $h\neq k$.~Then $M$ is a $\mathbb{G}$-basis.
\end{thm}
\begin{proof} Since we already proved that $\mathbb{G}$ is strongly represented by $M$ (see Proposition \ref{prop-strong G}), then it remains to show that its components are pairwise orthogonal.~In particular we verify  that for any pair of indices $h,k$ with $h\neq k$
 \begin{align*}
 [M^{h},M^{k}]\equiv 0,  \ \ \text{a.s.}
 \end{align*}
In fact by (\ref{eq-martingales}) and by linearity of the covariation operator
 \begin{align*}
  [M^h,M^k]=[N^{h},N^{k}]-[N^{h},N^{k,p}]-[N^{h,p},N^{k}]+[N^{h,p},N^{k,p}].
  \end{align*}
  Now observe that all the involved processes in the right hand side are bounded variation processes, and, as well known, for this kind of processes the quadratic covariation coincides with the sum of common jumps.\\
 Then, recalling that by construction $N^{h}$ and $N^{k}$ with $h\neq k$ do not jump together and applying Proposition \ref{prop-orth} and Corollary \ref{rem:1} to the processes $N^{h}$ and $N^{k}$,  we easily get $[M^h,M^k]\equiv 0$ a.s.
\end{proof}
   \noindent\begin{remark}\label{orth}
    Fixed $h$ let us denote by $\Big(\zeta^{h,dp,l}_{j}\Big)_j$ a sequence of $\mathbb{G}$-predictable random times enveloping the jump time indexed by $l$ of $M^{h,dp}$, the accessible part of $M^{h}$ (see Theorem 1-4 in \cite{yoeurp76}).~Analogously, fixed $k\neq h$, let us consider a sequence $\Big(\zeta^{k,dp,m}_{i}\Big)_i$
   of $\mathbb{G}$-predictable random times enveloping the jump time indexed by $m$ of $M^{k,dp}$, the accessible part of $M^{k}$.~Then the Assumption \ref{ass:NOPJT} for the processes $N^{h}$ and $N^{k}$ holds when
$$\Big(\bigcup_{l,j} [[\zeta^{h,dp,l}_{j}]]\Big) \cap    \Big(\bigcup_{m,i}[[\zeta^{k,dp,m}_{i}]]\Big)$$
is a set of null measure.
\end{remark}
\begin{remark}\label{confronto2 DT Jeanbl}
Let us consider again the case of two point processes $X$ and $H$ (see Remark \ref{confronto1 DT Jeanbl}).~Point (ii) of Theorem 3.5 in \cite{ditella-jean-21}  states that, denoting by $Z^1,Z^2,Z^3$ the $\mathbb{G}$-local martingales obtained by  $\mathbb{G}$-dual predictable compensation of $X-[X,H]$,  $H-[X,H]$ and $[X,H]$  respectively, then the triplet $(Z^1,Z^2,Z^3)$  enjoys the $\mathbb{G}$-PRP.~Moreover in the same paper, when $H$ coincides with the occurrence process of a finite random time $\tau$, that is $H=\mathds{1}_{[[\tau,+\infty[[}$, different conditions are considered making $(Z^1, Z^2, Z^3)$ a triplet of pairwise orthogonal local martingales, in other words a $\mathbb{G}$-basis.~It turns out that pairwise orthogonality holds in particular in two cases: the first one is when at least one among $X$ and $H$ is $\mathbb{G}$-quasi-left continuous, so that at least one of the  $\mathbb{G}$-predictable compensators is a continuous process; the second one is when $\tau$ avoids $\mathbb{F}$-stopping times (see Theorem 4.4. point (ii) and Corollary 4.5. in \cite{ditella-jean-21}).~Remark \ref{on NOPJT} explains that both situations are covered by our Assumption \ref{ass:NOPJT} and therefore by Theorem \ref{thm-ort-basis}.~In particular one has to recall that when $\tau$ avoids $\mathbb{F}$-stopping times, and therefore $\mathbb{F}$-predictable stopping times, then $\tau$ is $\mathbb{G}$-totally inaccessible (see page 65 in \cite{Jeulin}).\\Finally Theorem \ref{thm-ort-basis} answers very quickly  to the question posed by the first example in Section 4.3 of \cite{ditella-jean-21}.~The example deals with the representation in strong form of the progressive enlargement  of the reference filtration of an homogeneous Poisson process $X$ by the discrete random time $X_T+1$, where $T\in\mathbb{R}^+\setminus\{0\}$.~Again Assumption \ref{ass:NOPJT} holds true, since $X_T+1$ avoids the jump times of $X$ (see again Remark \ref{on NOPJT}).
\end{remark}
\section*{Acknowledgments}
\noindent The authors acknowledge the MIUR Excellence Department Project
awarded to the Department of Mathematics, University of Rome Tor
Vergata, CUP E83C18000100006.
\bibliography{biblio-spa}
\end{document}